\documentclass{amsart}
\usepackage[cp1251]{inputenc}
\usepackage{amsmath,amsopn,amsthm,amscd,amsfonts,longtable,amssymb,srcltx,comment}
\usepackage[final]{graphicx}
\usepackage{euscript}

\newtheorem{theorem}{Theorem}[section]
\newtheorem*{Bellows conjecture}{Bellows conjecture}
\newtheorem{corollary}{Corollary}[section]

\newtheorem{lemma}{Lemma}[section]
\newtheorem{proposition}{Proposition}[section]

\theoremstyle{definition}
\newtheorem{example}{Example}[section]
\newtheorem{remark}{Remark}[section]

\DeclareMathOperator{\per}{per}

\begin{document}

\title[Volumes of polytopes in spaces of constant curvature]{Volumes of polytopes \\ in spaces of constant curvature}
\author{Nikolay Abrosimov${}^1$}
\author{Alexander Mednykh${}^2$} 


\thanks{The work is supported by Russian Foundation for Basic Research (grants 12-01-00210, 12-01-33058 and 12-01-31006) and Council for Grants of President of the Russian Federation (projects MK-4447.2012.1 and SS-921.2012.1).}

\subjclass{Primary 51M20; Secondary 51M25, 51M09, 52B15}

\keywords{Volumes of polyhedra, constant curvature spaces, Sforza formula, Seidel problem, Brahmagupta formula}

\maketitle

\begin{center}
{\it ${}^{1,2}$ Sobolev Institute of Mathematics

${}^{1,2}$ Novosibirsk State University

${}^1$ abrosimov@math.nsc.ru, ${}^2$ mednykh@math.nsc.ru}
\end{center}

\begin{abstract}
We overview the volume calculations for polyhedra in Euclidean, spherical and hyperbolic spaces. We prove the Sforza formula for the volume of an arbitrary tetrahedron in $H^3$ and $S^3$. We also present some results, which provide a solution for Seidel problem on the volume of non-Euclidean tetrahedron. Finally, we consider a convex hyperbolic quadrilateral inscribed in a circle, horocycle or one branch of equidistant curve. This is a natural hyperbolic analog of the cyclic quadrilateral in the Euclidean plane. We find a few versions of the Brahmagupta formula for the area of such quadrilateral. We also present a formula for the area of a hyperbolic trapezoid.
\end{abstract}

{\normalsize \tableofcontents }

\section{Volumes of Euclidean polyhedra}

Calculating volumes of polyhedra is a classical problem, that has been well known since Euclid and remains relevant nowadays. This is partly due to the fact that the volume of a fundamental polyhedron
is one of the main invariants for a three-dimensional manifold.

One of the first results in this direction was obtained by Tartaglia (1499--1557), who 
had described an algorithm for calculating the height of a tetrahedron with some concrete lengths of its edges. The formula which expresses the volume of an Euclidean tetrahedron in terms of its edge lengths was given by Euler (see \cite{Usp}, p.~256). The multidimensional analogue of this result is known as the Cayley--Menger formula (see \cite{Sommer}, p.~124).

\begin{theorem} [Tartaglia, XVI AD] Let
$T$ be an Euclidean tetrahedron with edge lengths
$d_{ij}$, $1\leq i<j\leq4$. Then the volume $V=V(T)$ is given by the formula
$$
288\,V^2=\left|%
\begin{array}{ccccc}
  0 & 1 & 1 & 1 & 1 \\
  1 & 0 & d_{12}^2 & d_{13}^2 & d_{14}^2 \\
  1 & d_{21}^2 & 0 & d_{23}^2 & d_{24}^2 \\
  1 & d_{31}^2 & d_{32}^2 & 0 & d_{34}^2 \\
  1 & d_{41}^2 & d_{42}^2 & d_{43}^2 & 0 \\
\end{array}%
\right|.$$
\end{theorem}

In the above relation the volume is evaluated as a root of the quadratic equation whose coefficients are polynomials with integer coefficients in edge lengths. Surprisingly, but this result can be generalized to an arbitrary Euclidean polyhedron. About fifteen years ago, I.~Kh.~Sabitov \cite{Sabitov195} proved the corresponding theorem for any polyhedron in $R^3$.
Then Robert Connelly, Idzhad Sabitov and Anke Walz gave the second proof for general orientable $2$-dimensional polyhedral surfaces using the theory of places instead of resultants \cite{CSW}.

\begin{theorem} [Sabitov, 1996; Connelly, Sabitov, Walz, 1997]\label{Sabitov}
Let $P$ be Euclidean polyhedron with triangular faces and edge lengths $d_{ij}$. Then the volume
$V(P)$ is a root of a polynomial whose coefficients depend on $d_{ij}^2$ and
combinatorial type of $P$ only.
\end{theorem}


Note that the explicit form of the above mentioned polynomial is known only in some special cases, in particular for octahedra with symmetries \cite{GMS}. On the other hand, Theorem~\ref{Sabitov} allows to prove the well-known Bellows conjecture stated by R.~Connelly, N.~Kuiper and D.~Sullivan \cite{Con80}.

\begin{Bellows conjecture} [Connelly, Kuiper, Sullivan, late 1970s]
The generalized volume of a flexible polyhedron does not change when it is bending.
\end{Bellows conjecture}

Recall \cite{Sabitov195} that the \emph{generalized volume} of an oriented geometrical polyhedron is equal to the sum of volumes of consistently oriented tetrahedra with a common vertex and with bases on the faces of the polyhedron. \emph{Bending} of a polyhedron is a continuous isometric deformation, provided the rigidity of the faces.


Cauchy rigidity theorem \cite{Cauchy} states that a convex polyhedron with rigid faces is rigid itself. For non-convex polyhedra it is not true; there are examples of flexible polyhedra among them. The first example of a flexible polyhedron was constructed by Bricard \cite{Bricard}. It was a self-intersecting octahedron. The first example of a flexible polyhedron embedded into Euclidean $3$-space was presented by Connely \cite{Con78}. The smallest example of such polyhedron is given by Steffen. It has $14$ triangular faces and $9$ vertices (see Fig.~1). 

\begin{figure}[htbp]
\begin{center}
\includegraphics[scale=0.5]{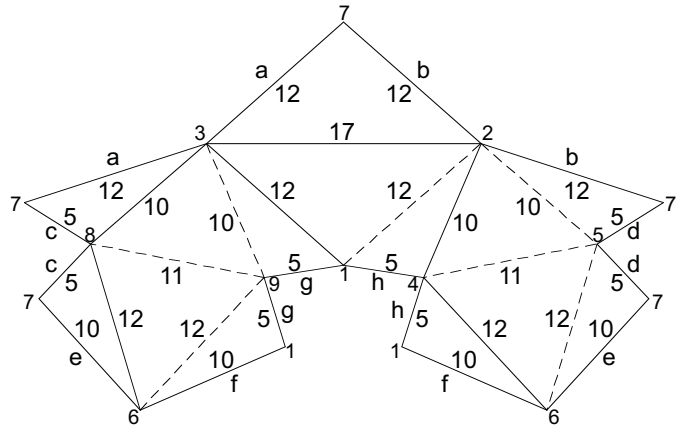}
\end{center}
\caption{Development of Steffen's flexible polyhedron}
\end{figure}


A flexible polyhedron keeps the same combinatorial type and the same set of edge lengths under bending.
Then, by Theorem~\ref{Sabitov} the volume of this polyhedron can take only a finite number of values corresponding to the roots of a polynomial. Since the bending of a polyhedron is a continuous isometric deformation, the volume is constant.


A few months ago a new result by A.~A.~Gaifullin \cite{Gaif} was published in arXiv. He proved an analog of Theorem~\ref{Sabitov} for the generalized volume of a $4$-dimensional polyhedron. In the spherical space the analog of Theorem~\ref{Sabitov} is not true (see \cite{Alex97}) while in the hyperbolic space the question is still open. One will see below that in many cases the volume of a non-Euclidean polyhedron is not expressed in terms of elementary functions.


\section{Volumes of non-Euclidean tetrahedra}

In the hyperbolic and spherical cases, the situation becomes more complicated. Gauss, who was one of the creators of hyperbolic geometry, referred to the calculation of volumes in non-Euclidean geometry as {\it die Dschungel}.

It is well known that a tetrahedron in $S^3$ or $H^3$ is determined by an ordered set of its dihedral angles up to isometry of the space (see, e.g., \cite{Avs}). Recall that in the Euclidean case this is true up to simility.

\subsection{Orthoschemes in $S^3$ and $H^3$}

Volume formulas for non-Euclidean tetrahedra in some special cases has been known
since Lobachevsky, Bolyai and Schl\"{a}fli. For example, Schl\"{a}fli \cite{Schla} found the volume of an orthoscheme in $S^3$. Recall that an \emph{orthoscheme} is an $n$-dimensional simplex defined by a sequence of edges $(v_0,v_1), (v_1,v_2),\ldots,(v_{n-1},v_n)$ that are mutually orthogonal. In three dimensions, an orthoscheme is also called a \emph{birectangular tetrahedron} (see Fig.~2).  

\begin{figure}[htbp]
\begin{center}
\includegraphics[scale=0.45]{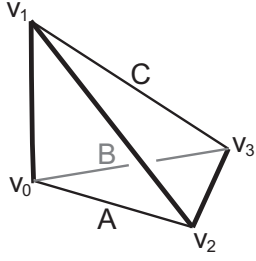}
\end{center}
\caption{Orthoscheme $T=T(A,B,C)$ with essential dihedral angles $A,B,C$, all other angles are $\displaystyle\frac{\pi}{2}$}
\end{figure}

\begin{theorem} [Schl\"{a}fli, 1858]\label{SchOrt}
Let $T$ be a spherical orthoscheme with essential dihedral angles $A, B, C$. Then the volume $V=V(T)$ is given by the formula  
\begin{equation*}
\begin{split}
&V=\displaystyle\frac{1}{4}S(A, B, C),  \,\,\,where \,\,\,S\left(\frac{\pi}{2}-x,y,\frac{\pi}{2}-z\right)=\widehat{S}(x,y,z)=\\
&\sum\limits_{m=1}^{\infty}\left(\frac{D-\sin x\sin z}{D+\sin x\sin z}\right)^m\frac{\cos 2mx-\cos 2my+\cos 2mz-1}{m^2}
-x^2+y^2-z^2,
\end{split}
\end{equation*}
where $D\equiv\sqrt{\cos^2x\cos^2z-\cos^2y}.$
\end{theorem}

The function $S(x,y,z)$ appeared in Theorem~\ref{SchOrt} is called the \emph{Schl\"{a}fli function}. The volume of a hyperbolic orthoscheme was obtained independently by Janos Bolyai \cite{Boly} and Nikolai Lobachevsky \cite{Loba}. The following theorem represents a result of Lobachevsky as a quite simple formula. In such form it was given by G.~S.~M.~Coxeter \cite{Coxe}.


\begin{theorem} [Lobachevsky, 1835; Coxeter, 1935]
Let $T$ be a hyperbolic orthoscheme with dihedral angles $A,B,C$. Then the volume $V=V(T)$ is given by the formula $V=\displaystyle\frac{i}{4}S(A, B, C)$, where $S(A, B, C)$ is the Schl\"{a}fli function.
\end{theorem}

J.~Bolyai (see, e.g., \cite{Boly}) found the volume of a hyperbolic orthoscheme in terms of planar angles and an edge length. Consider a hyperbolic orthoscheme $T$ with essential dihedral angles along $AC, CB, BD$, and all other dihedral angles equaled to $\displaystyle\frac{\pi}{2}$ (Fig.~3).

\begin{figure}[htbp]
\begin{center}
\includegraphics[scale=0.3]{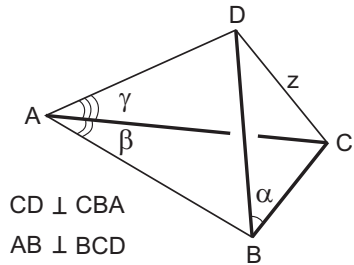}
\end{center}
\caption{Orthoscheme $T$ given by planar angles $\alpha,\beta,\gamma$ and edge length $z$}
\end{figure}

\begin{theorem} [Bolyai, 1832]
Let $T$ be a hyperbolic orthoscheme with planar angles $\alpha,\beta,\gamma$ and edge length $z$. Then the volume $V=V(T)$ is given by the formula 
$$
V=\frac{\tan\gamma}{2\tan\beta}\int\limits_0^z\frac{u\,\sinh u\, du}{\left(\frac{\cosh^2 u}{\cos^2\alpha}-1\right)\sqrt{\frac{\cosh^2 u}{\cos^2\gamma}-1}}.
$$
\end{theorem}


\subsection{Ideal tetrahedra}

An \emph{ideal tetrahedron} is a tetrahedron with all vertices at infinity. Opposite dihedral angles of an ideal tetrahedron are pairwise equal and the sum of dihedral angles at the edges adjacent to one vertex is $A+B+C=\pi$ (see Fig.~4).


\begin{figure}[htbp]
\begin{center}
\includegraphics[scale=0.3]{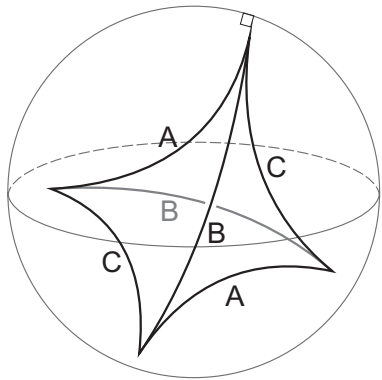}
\end{center}
\caption{Ideal tetrahedron $T=T(A,B,C)$ with dihedral angles $A,B,C$}
\end{figure} 


The volume of an ideal tetrahedron has been known since Lobachevsky \cite{Loba}. J.~Milnor presented it in very elegant form \cite{Milnor}.


\begin{theorem} [Lobachevsky, 1835; Milnor, 1982]\label{Miln}
Let $T$ be an ideal hyperbolic tetrahedron with dihedral angles $A, B$ and $C$. Then the volume $V=V(T)$ is given by the formula
$$V=\Lambda(A)+\Lambda(B)+\Lambda(C),$$
where $\Lambda(x)=-\displaystyle\int\limits_0^x \log|2\sin t|\,dt$ is the Lobachevsky function.
\end{theorem}

More general case of tetrahedron with at least one vertex at infinity was investigated by E.~B.~Vinberg \cite{Avs}.

\subsection{General hyperbolic tetrahedron}

Despite the fact that partial results on the volume of non-Euclidean tetrahedra were known a volume formula for hyperbolic tetrahedra of general type remained unknown until recently. A general algorithm for obtaining such a formula was indicated by W.--Yi.~Hsiang in \cite{Hsiang}. A complete solution was obtained more than ten years later. In paper by Korean mathematicians Y.~Cho and H.~Kim \cite{ChoKim} proposed a general formula. However, it was asymmetric with respect to permutation of angles. The next advance was achieved by Japanese mathematicians. First, J.~Murakami and Y.~Yano \cite{MurYan} proposed a formula expressing the volume by dihedral angles in symmetric way. A year later, A.~Ushijima \cite{Ush} presented a simple proof of the Murakami--Yano's formula. He also investigated the case of a truncated hyperbolic tetrahedron.

It should be noted that in all these studies the volume is expressed as a linear combination of 16 dilogarithms or Lobachevsky functions. The arguments of these functions depend on the dihedral angles of the tetrahedron, and some additional parameter, which is the root of some quadratic equation with complex coefficients in a sophisticated form.

The geometric meaning of this formula was explained by P.~Doyle and G.~Leibon \cite{Leib} in terms of Regge symmetry. A clear description of these ideas and a complete geometric proof of this formula was given by Yana Mohanty \cite{MohPhD}. In particular, she was able to prove the equivalence of Regge symmetry and homogeneity (scissors congruence) in the hyperbolic space \cite{Moh2}.

In 2005, D.~A.~Derevnin and A.~D.~Mednykh \cite{DM2005} proposed the following integral formula for the volume of a hyperbolic tetrahedron.

\begin{theorem} [Derevnin, Mednykh, 2005]
Let $T(A,B,C,D,E,F)$ be a compact hyperbolic tetrahedron with dihedral angles $A,B,C,D,E,F$. Then the volume $V=V(T)$ is given by the formula
$$
V=-\frac14
\int\limits_{z_1}^{z_2}\log\frac{\cos{\frac{A+B+C+z}{2}}\cos{\frac{A+E+F+z}{2}}
\cos{\frac{B+D+F+z}{2}}\cos{\frac{C+D+E+z}{2}}}
{\sin{\frac{A+B+D+E+z}{2}}\sin{\frac{A+C+D+F+z}{2}}\sin{\frac{B+C+E+F+z}{2}}\sin{\frac{
z}{2}}}\,dz\,,
$$
where $z_1$ and $z_2$ are the roots of the integrand satisfied the condition $0< z_2 - z_1 < \pi\,.$ More precisely,
\begin{align*}
z_1&=\arctan\displaystyle\frac{k_3}{k_4}-\arctan\frac{k_1}{k_2}\,,\,\,\,\,\, z_2=\pi-\arctan\frac{k_3}{k_4}-\arctan\frac{k_1}{k_2}\,,\,\,\,\,\,\ where
\\
k_1&=-\cos S-\cos (A+D)-\cos (B+E)-\cos (C+F)-\cos (D+E+F)-
\\
&-\cos(D+B+C)-\cos (A+E+C)-\cos (A+B+F)\,,
\\
k_2&=\sin S+\sin (A+D)+\sin (B+E)+\sin (C+F)+\sin (D+E+F)+
\\
&+\sin (D+B+C)+\sin (A+E+C)+\sin (A+B+F)\,,
\\
k_3&=2\,(\sin A\sin D+\sin B\sin E+\sin C\sin F)\,,
\\
k_4&=\sqrt{k_1^2+k_2^2-k_3^2}\,,
\\
S&=A+B+C+D+E+F\,.
\end{align*}
\end{theorem}

\subsection{Sforza formula for non-Euclidean tetrahedron}

Surprisingly, but more than 100 years ago, in 1907, Italian mathematician Gaetano Sforza (or Scorza in some of his papers) found a fairly simple formula for the volume of a non-Euclidean tetrahedron. This fact has became widely known after a discussion between the second named author and J.~M.~Montesinos at the conference in El Burgo de Osma (Spain), August, 2006. Unfortunately, the outstanding work of Sforza \cite{Sforza} published in Italian has been completely forgotten. 

The original arguments by Sforza are based on some identity given by H.~W.~Richmond \cite{Rich} besides of Schl\"{a}fli formula. Then he used a Pascal's equation \cite{Pasc} for minors of Gram matrix and some routine calculations. In this section we provide a new proof of Sforza formula for the volume of an arbitrary tetrahedron in $H^3$ or $S^3$. The idea of this proof belongs to authors, it was never published before. 

Consider a hyperbolic (or a spherical) tetrahedron $T$ with dihedral angles $A,B,C,$ $D,E,F$. Assume that $A, B, C$ are dihedral angles at the edges adjacent to one vertex and $D, E, F$ are opposite to them correspondingly (see Fig.~5). Then the Gram matrix $G(T)$ is defined as follows

\begin{figure}[htbp]
\begin{center}
\includegraphics[scale=0.45]{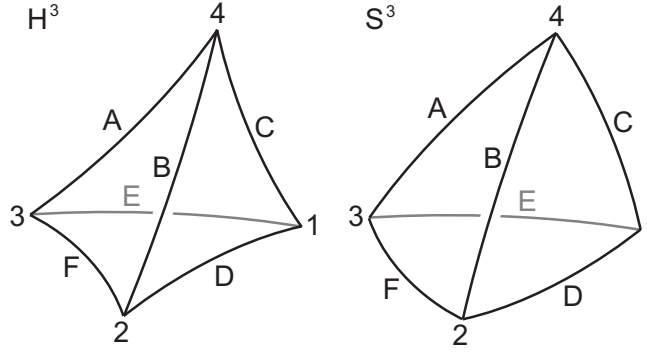}
\end{center}
\caption{Tetrahedron $T=T(A,B,C,D,E,F)$ in $H^3$ or $S^3$}
\end{figure}

$$G = \left(%
\begin{array}{cccc}
  1 & - \cos A & - \cos B & - \cos F \\
  - \cos A & 1 & - \cos C & - \cos E \\
  - \cos B & - \cos C & 1 & - \cos D \\
  - \cos F & - \cos E & - \cos D & 1 \\
\end{array}%
\right).$$

Denote by  $\mathcal{C}=\langle c_{ij} \rangle _{i,j=1,2,3,4}$ a
matrix of cofactors $c_{ij}=(-1)^{i+j} G_{ij}$,
where $G_{ij}$ is $(i,j)$-th minor of $G$.

In the following proposition we collect some known results about hyperbolic tetrahedron (see, e.g., \cite{Ush}).

\begin{proposition}\label{Ushi}
Let $T$ be a compact hyperbolic tetrahedron. Then
\begin{eqnarray}\label{eqo}
\text{\upshape(i)} &  & \det G<0; \nonumber
\\ \text{\upshape(ii)} &
& c_{ii}>0,\,i\in\{1,2,3,4\}; \nonumber
\\ \text{\upshape(iii)} &
& c_{ij}>0,\,i\ne j;\,i,j\in\{1,2,3,4\}; \nonumber
\\ \text{\upshape(iv)} & & \cosh{\ell_{ij}}=\frac{c_{ij}}{\sqrt{c_{ii}
c_{jj}}},\nonumber
\end{eqnarray}
where $\ell_{ij}$ is a hyperbolic length of the edge joining vertices
$i$ and $j$.
\end{proposition}

Further, we need the following assertion due to Jacobi (\cite{Pras94},
Theorem 2.5.1, p.~12).

\begin{proposition}[Jacobi theorem] \label{Jaco}
Let $M=\langle m_{ij} \rangle _{i,j=1,...,n}$ be a matrix and
$\Delta= \det M$ is determinant of $M$. Denote by
$\mathcal{C}=\langle c_{ij} \rangle _{i,j=1,...,n}$ the matrix of cofactors $c_{ij}=(-1)^{i+j}\det M_{ij}$, where $M_{ij}$ is
$(n-1)\times(n-1)$ minor obtained by removing $i$-th row and
$j$-th column of the matrix $M$. Then for any $k$, $1\leq k\leq
n-1$ we have
\begin{equation*}
 \det \langle c_{ij} \rangle
_{i,j=1,...,k}=\Delta^{k-1}\det\langle m_{ij} \rangle
_{i,j=k+1,...,n}\,.
\end{equation*}
\end{proposition}

One of the main tools for volume calculations in $H^3$ and $S^3$ is the classical Schl\"{a}fli formula. 

\begin{proposition}[Schl\"{a}fli formula]\label{SchFor}
Let $X^n$ be a space of constant curvature $K$. Consider a family of convex polyhedra $P$ in $X^n$ depending on one or more parameters in a differential manner and keeping the same combinatorial type. Then the differential of the volume $V=V(P)$ satisfies the equation
$$(n-1)K\,dV =\sum\limits_{F}V_{n-2}(F)\,d\theta(F),$$
where the sum is taken over all $(n-2)$-facets of $P$, $V_{n-2}(F )$ is $(n-2)$-dimensional volume of $F$, and $\theta(F)$ is the interior angle along $F$.
\end{proposition}

In the classical paper by Schl\"{a}fli \cite{Sch1858-60} this formula was proved for the case of a spherical $n$-simplex. For the hyperbolic case, it was obtained by H.~Kneser \cite{Kn1936} (for more details, see also \cite{Avs} and \cite{M94}). In the Euclidean case this formula reduces to the identity $0=0$.

Now we are able to prove the following theorem.

\begin{theorem}[Sforza formula in $H^3$] 
Let $T$ be a compact hyperbolic tetrahedron with Gram matrix $G$. Consider $G=G(A)$ as a function of dihedral angle $A$. Then the volume $V=V(T)$ is given by the formula

$$V=\frac14\,\int\limits_{A_0}^A\log\frac{c_{34}(A)-\sqrt{-\det\,G(A)} \sin A}{c_{34}(A)+
\sqrt{- \det\,G(A)} \sin A}\,d A,$$
where $A_0$ is a suitable root of the equation $\det\,G(A)=0$ and $c_{34}=c_{34}(A)$ is $(3,4)$-cofactor of the matrix $G(A)$.
\end{theorem}
{\sc Proof} (Abrosimov, Mednykh){\bf.} 
Let $\Delta =\det\,G.$  By the Jacobi theorem (Proposition~\ref{Jaco})
applied to Gram matrix $G$ for $n=4$ and $k=2$ we obtain
$$c_{33}c_{44}-c_{34}^2=\Delta (1-\cos^2 A). $$ By the Cosine rule (Proposition~\ref{Ushi}, (iv)) we get
$ \displaystyle\cosh \ell_A=\frac{c_{34}}{\sqrt{c_{33}c_{44}}}.$ Hence,

$$ \sinh \ell_A=\sqrt{\cosh^2 \ell_A-1}=\sqrt{\frac{c_{34}^2-c_{33}c_{44}}{c_{33}c_{44}}}
 =\frac{\sqrt{ -\Delta}\,\sin A} {\sqrt{c_{33}c_{44}}}.$$
\noindent Since $\exp(\pm \ell_A)=\cosh \ell_A \pm \sinh \ell_A$
we have
$$\exp(\pm  \ell_A)=\frac{c_{34}\pm \sqrt{-\Delta}\,\sin A}{\sqrt{c_{33}c_{44}}} .$$
Hence,
$$\exp(2\,\ell_A)=\frac{\exp(\ell_A)}{\exp(-\ell_A)}= \frac{c_{34}+ \sqrt{-\Delta}\,\sin A}{c_{34}- \sqrt{-\Delta}\,\sin A}$$ and
$$  \ell_A = \frac12\,\log \frac{c_{34}+ \sqrt{-\Delta}\,\sin A}{c_{34}- \sqrt{-\Delta}\,\sin A}.$$
By the Schl\"{a}fli formula (Proposition~\ref{SchFor}) for $V=V(T)$ we have
$$-dV=\frac12 \sum \limits_{\alpha} \ell_{\alpha}\,d\alpha, \,\,\,\alpha
\in \{A, B, C, D, E, F\}.$$ 
By assumption that angle $A$ is variable and all other angles are constant, we get
$$-dV=\frac12\, \ell_A\,dA.$$ 
Note that $\det G\rightarrow0$ as $A\rightarrow A_0$. Thus, $V\rightarrow0$ as $A\rightarrow A_0$.
Then integrating both sides of the equation we obtain
$$V=\int\limits_{A_0}^A\left(-\frac{\ell_A}{2}\right)\,d A=
\frac14\,\int\limits_{A_0}^A\log\frac{c_{34}-\sqrt{-\det G}\,\sin A}{c_{34}+\sqrt{-
\det G}\,\sin A}\,d A,$$
where the lower limit $A_0$ is a suitable root of the equation
$\det G(A)=0$.\,\,\,$\Box$

In the following proposition we collect some known results about spherical tetrahedron (see, e.g., \cite{Luo}).

\begin{proposition}\label{Luo}
Let $T$ be a spherical tetrahedron. Then
\begin{eqnarray}\label{eqo}
\text{\upshape(i)} &  & \det G>0; \nonumber
\\ \text{\upshape(ii)} &
& c_{ii}>0, i=1,2,3,4; \nonumber
\\ \text{\upshape(iii)} & & \cos{\ell_{ij}}=\frac{c_{ij}}{\sqrt{c_{ii}
c_{jj}}},\nonumber
\end{eqnarray}
where $\ell_{ij}$ is a length in $S^3$ of the edge joining vertices
$i$ and $j$.
\end{proposition}

The next theorem presents a spherical version of the Sforza formula.

\begin{theorem}[Sforza formula in $S^3$] 
Let $T$ be a spherical tetrahedron with Gram matrix $G$. Consider $G=G(A)$ as a function of dihedral angle $A$. Then the volume $V=V(T)$ is given by the formula

$$V=\frac{1}{4i}\,\int\limits_{A_0}^A\log\frac{c_{34}(A)+i\sqrt{\det\,G(A)} \sin A}{c_{34}(A)-
i\sqrt{\det\,G(A)} \sin A}\,dA,$$
where $A_0$ is a suitable root of the equation $\det\,G(A)=0$ and $c_{34}=c_{34}(A)$ is $(3,4)$-cofactor of the matrix $G(A)$.
\end{theorem}

The proof is similar to one given in the hyperbolic case.

\section{Seidel conjecture}

In 1986, J.~J.~Seidel \cite{Seidel} conjectured that the volume of an ideal hyperbolic tetrahedron can be expressed as a function of the determinant and the permanent of its Gram matrix. Recall that the formula expressing the volume of such tetrahedron in terms of dihedral angles has been known since Lobachevsky and Bolyai (Theorem~\ref{Miln}). In spite of this, the Seidel problem had not been solved for a long time. 10 years later, a strengthened version of Seidel conjecture was suggested by Igor Rivin and Feng Luo. They supposed that the volume of a non-Euclidean tetrahedron (hyperbolic or spherical one) depends only on the determinant of its Gram matrix. It was shown in \cite{A09} and \cite{A10} that the strengthened conjecture is false, while Seidel conjecture is true within certain conditions.

Consider a non-Euclidean tetrahedron $T$ with dihedral angles $A,B,C,D,E,F$ in $S^3$ or $H^3$ (Fig.~5). Denote vertices by numbers $1,2,3,4$. Let $A_{ij}$ denote a dihedral angle at the edge joining vertices $i$ and $j$. For convenience, we set $A_{ii} =\pi$ for $i=1,2,3,4$. It is well known [4] that, in the hyperbolic and spherical spaces, the tetrahedron $T$ is uniquely (up to isometry) determined by its Gram matrix 
$$G = \langle-\cos A_{ij}\rangle_{i,j=1,2,3,4}
= \left(%
\begin{array}{cccc}
  1 & - \cos A & - \cos B & - \cos F \\
  - \cos A & 1 & - \cos C & - \cos E \\
  - \cos B & - \cos C & 1 & - \cos D \\
  - \cos F & - \cos E & - \cos D & 1 \\
\end{array}%
\right).$$


Recall that the permanent of a matrix $M=\langle m_{ij}\rangle_{ij=1,2,\dots,n}$ is defined by 
$${\rm per}M=\sum\limits_{i=1}^n m_{ij}\,{\rm per}M_{ij},\,\,\,\,{\rm per}(m_{ij})=m_{ij},$$
where $M_{ij}$ is the matrix obtained from $M$ by removing $i$-th row and the $j$-th column.
Conditions for the existence of spherical and hyperbolic tetrahedra in terms of Gram matrices are given in \cite{Luo} and \cite{Ush}, respectively.

\subsection{Strengthened conjecture}
In this section we provide a counterexample to Seidel strengthened problem.
In the spherical case, the answer is given by the following theorem \cite{A09}.

\begin{theorem}[Abrosimov, 2009]
There exists a one-parameter family of spherical tetrahedra with unequal volumes and the same determinant of Gram matrix.
\end{theorem}

We prove this theorem by constructing such a family. Consider a tetrahedron $T(A, D)$ with two opposite dihedral angles equaled to $A$ and $D$, and the remaining dihedral angles equaled to $\displaystyle\frac{\pi}{2}$ (Fig.~6). 

\begin{figure}[htbp]
\begin{center}
\includegraphics[scale=0.25]{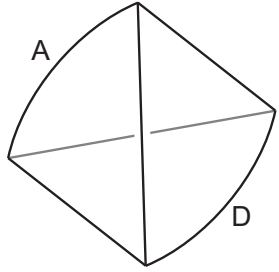}
\end{center}
\caption{Tetrahedron $T=T(A,D)$}
\end{figure}

It is easy to show that the volume of such a tetrahedron equals $\displaystyle\frac{A\,D}{2}$ and
the determinant of its Gram matrix is $\det G=\sin^2 A\,\sin^2\,D$. Among all tetrahedra $T(A,D)$ with $0<A,D<\pi$, we choose a family of tetrahedra $$\displaystyle T_c(A,D)=T\left(A,\arcsin\frac{c}{\sin A}\right)$$ whose Gram matrices have the same determinant $\det G=c^2$, where $c$ is a constant satisfying the inequalities $0<c<\min\{\sin A,\,\sin D\}$.

The volume of such tetrahedra is given by equation
$$V(T_c)=\frac{A}{2}\,\arcsin\frac{c}{\sin A}.$$
Thus, it depends not only on the constant $c$, but also on the value of the free parameter $A$.

We have constructed the required family of tetrahedra.
In the hyperbolic case, we have failed to construct an elementary counterexample to Seidel strengthened conjecture. Nevertheless, a similar theorem is also valid \cite{A10}.

\begin{theorem}[Abrosimov, 2009]
There exists a one-parameter family of hyperbolic tetrahedra with unequal volumes and the same determinant of Gram matrix.
\end{theorem}

The proof of this theorem is based on the following considerations. Consider any hyperbolic tetrahedron $T$ with dihedral angles $A,B,C,D,E,F$. Assume that the angles $A,B,C$ are along the edges adjacent to one vertex, and $D,E,F$ are opposite to them. We fix all dihedral angles except two opposite ones, say $A$ and $D$. Since the set of hyperbolic tetrahedra is open (see \cite{Ush}, \cite{Luo}), it follows that, varying $A$ and $D$ within sufficiently small limits, we obtain hyperbolic tetrahedra. In the set $T(A,D)$ of such tetrahedra, we choose a family of tetrahedra $T_c(A,D)$ with the same Gram determinant $\det G=–c^2<0$. The latter condition means that the differential of the function $\det G$ is zero. Since the angles $A$ and $D$ vary and the remaining angles are fixed, it follows that
$$-d\det G=2c_{12}\sin A\,dA+2c_{34}\sin D\,dD=0,$$
where $c_{ij}$ is $(i,j)$-th cofactor of the matrix $G$. Due to this relation, we can treat the angle $D$ as a function of the angle $A$. We have
$$\frac{dD}{dA}=-\frac{c_{12}\sin A}{c_{34}\sin D}.$$

The derivative of the volume as a composite function of the angle $A$ equals
$$\frac{dV}{dA}=\frac{\partial V}{\partial A}+\frac{\partial V}{\partial D}\,\frac{dD}{dA}.$$

According to classical Schl\"{a}fli formula (Proposition~\ref{SchFor}), we have
$$\frac{\partial V}{\partial A}=-\frac{\ell_A}{2},\,\,\,\,\frac{\partial V}{\partial D}=-\frac{\ell_D}{2},$$
where $\ell_A$ and $\ell_D$ are the lengths of the corresponding edges of the tetrahedron.

In turn, the lengths of the edges can be expressed in terms of dihedral angles (see \cite{Ush}, \cite{MP06})
\begin{align*}
\ell_A&={\rm arccoth}\frac{\sqrt{-\det G}\,\sin A}{c_{34}},\\
\ell_D&={\rm arccoth\frac{\sqrt{-\det G}\,\sin D}{c_{12}}}.
\end{align*}

Comparing these expressions and performing simple calculations, we obtain
$$\frac{dV}{dA}=-\frac{\tanh\ell_A}{2}\left(\frac{\ell_A}{\tanh\ell_A}-\frac{\ell_D}{\tanh\ell_D}\right).$$

It is required in the theorem that the volume is varying according to parameter $A$, i.e., $\displaystyle\frac{dV}{dA}\ne0$, which is equivalent to condition $\ell_A\ne\ell_D$.

Thus, two sufficiently close tetrahedra from the family $T_c(A,D)$ have the same Gram determinant and unequal volume provided that $\ell_A\ne\ell_D$. It is not hard to construct an infinite family of tetrahedra satisfying the last condition for $A\ne D$. For example, this condition is satisfied by ``almost symmetric'' tetrahedra with angles $A\ne D, B=E,$ and $C=F$. Recall that, for fixed $c$, the family $T_c(A,D)$ still depends on one free parameter.

\subsection{Solution of Seidel conjecture}
The solution of Seidel problem, which was posed in \cite{Seidel}, is given by the following theorem \cite{A10}.

\begin{theorem}[Abrosimov, 2010]\label{Seidel2010}
In each of classes of acute or obtuse tetrahedra the volumes of ideal tetrahedra are determined uniquely by the determinant and permanent of Gram matrix.\footnote{ By an \emph{obtuse tetrahedron} we mean a tetrahedron with at least one dihedral angle $\displaystyle >\frac{\pi}{2}$.}
\end{theorem}

As is known (see, e.g., \cite{Milnor}), opposite dihedral angles of an ideal tetrahedron are pairwise equal and the sum of dihedral angles at the edges adjacent to one vertex is $A+B+C=\pi$ (Fig.~4).

Thus, we can take $C=\pi-A-B$. The Gram matrix has the form
$$G = \left(%
\begin{array}{cccc}
  1 & - \cos A & - \cos B & \cos (A+B) \\
  - \cos A & 1 & \cos (A+B) & - \cos B \\
  - \cos B &  \cos (A+B) & 1 & - \cos A \\
  \cos (A+B) & - \cos B & - \cos A & 1 \\
\end{array}%
\right).$$

We have
\begin{align*}
\det G&=-4\sin^2A\,\sin^2B\,\sin^2(A+B),\\
{\rm per}\,G&=4+4\cos^2A\,\cos^2B\,\cos^2(A+B).
\end{align*}

To prove Theorem~\ref{Seidel2010}, we show that the dihedral angles of an ideal tetrahedron are uniquely (up to a permutation) determined by $\det G$ and ${\rm per}\,G$ in each of the cases: for acute angled tetrahedron and for obtuse angled tetrahedron.

Without loss of generality, we can assume that $0<A\leq B\leq C=\pi-A-B$. Then the dihedral angles $A,B$ are a priori acute, and the angle $C$ is either acute or obtuse. In the former case, the tetrahedron under consideration is acute-angled, and in the latter case, it is obtuse-angled.

Let us introduce the new variables
$$x=\sin A\,\sin B,\,\,\,\,y=\cos A\,\cos B$$
and show that, for a fixed left-hand side, the solutions of the system of equations
\begin{align*}
-\frac{1}{4}\det G&=x^2(1-(y-x)^2),\\
\frac{1}{4}\,{\rm per}\,G&=y^2(y-x)^2+1
\end{align*}
correspond to one tetrahedron (up to isometry) in each of the two cases mentioned above.

Suppose that the system has a pair of different solutions $(a,b)$ and $(x,y)$. Then we have
\begin{align*}
a^2(1-(b-a)^2)&=x^2(1-(y-x)^2),\\
b^2(b-a)^2&=y^2(y-x)^2.
\end{align*}

The angle $C$ being acute means that $\cos A\,\cos B-\sin A\,\sin B = –\cos C<0$, i.e., both solutions satisfy the inequalities $b(b-a)<0$ and $x(x-y)<0$. If the angle $C$ is obtuse, then the reverse inequality. This allows us to take a square root in the second equation without loosing solutions:
\begin{align*}
a^2(1-(b-a)^2)&=x^2(1-(y-x)^2),\\
b\,(b-a)&=y\,(y-x).
\end{align*}

Expressing $x$ from the second equation and substituting the resulting expression into the first one, we obtain a sixth-degree polynomial equation in $y$. Fortunately, it decomposes into the linear factors and the biquadratic polynomial 
\begin{equation*}
\begin{split}
&(b-y)\cdot(b+y)\cdot\big(y^4-(a^2+a^4+2\,ab-2\,a^3b-b^2+2\,ab^3-b^4)\,y^2+\\
&a^4b^2-4\,a^3b^3+6\,a^2b^4-4\,ab^5+b^6\big)=0.
\end{split}
\end{equation*}


Thus, all solutions can be we found in radicals. Substituting the expressions for $x$ and $y$ in terms of dihedral angles, we see that different solutions of the system correspond to the same ideal tetrahedron $T(A,B,C)$ up to reordering the dihedral angles.

Note that, in Theorem~\ref{Seidel2010}, the assumption that the tetrahedron is either acute-angled or obtuse-angled cannot be dispensed with. This is demonstrated by the following example.

\begin{example}
Consider a pair of ideal tetrahedra, the obtuse-angled tetrahedron $T_1(s,s,\pi-2s)$ and the acute-angled tetrahedron $\displaystyle T_2\left(t,\frac{\pi-t}{2},\frac{\pi-t}{2}\right)$, where
\begin{align*}
s&=\arccos\frac{\sqrt{2+\sqrt{4+\sqrt{170\sqrt{17}-698}}}}{2\sqrt{2}},\\
t&=\arccos\frac{-1+\sqrt{17}+\sqrt{-26+10\sqrt{17}}}{8}.
\end{align*}

The determinants and the permanents of the Gram matrices of these tetrahedra coincide; they are
\begin{align*}
\det G(T_1)=\det G(T_2)&=\frac{107-51\sqrt{17}}{128},\\
\per G(T_1)=\per G(T_2)&=\frac{163+85\sqrt{17}}{128}.
\end{align*}

At the same time, the volumes of the tetrahedra $T_1$ and $T_2$ are different and equal $0.847365$ and $1.01483$, respectively.
\end{example}

\section{Heron and Brahmagupta formulas}

Heron of Alexandria (c. 60 BC) is credited with the following formula that relates the area $S$ of a triangle to its side lengths $a, b,c$ $$S^2 = (s - a)(s - b)(s - c)s,$$
where $s = (a + b + c)/2$ is the semiperimeter. For polygons with more than three sides, the side lengths do not in general determine the area, but they do if the polygon is convex and cyclic (inscribed in a circle). Brahmagupta, in
the seventh century, gave the analogous formula for a convex cyclic quadrilateral with side lengths $a,b,c,d$
$$S^2 = (s - a)(s - b)(s - c)(s - d),$$
where $s = (a + b + c + d)/2. $ See \cite{CoxCrei} for the elementary proof. An interesting consideration of the problem can be found in the M\"obius paper \cite{Moeb}. Independently, D.~P.~Robbins \cite{Robb} and V.~V.~Varfolomev \cite{Varf} found a way to generalize these formulas. The main idea of both papers was to determine the squared area $S^2$ as a root of an algebraic equation whose coefficients are integer polynomials in the squares of the side lengths. See also \cite{FedPak}, \cite{Conn} and \cite{Sab} for more detailed consideration.


In the present section of our work we deal with the hyperbolic plane instead of the Euclidean one. The hyperbolic plane under consideration is equipped by a Riemannian metric of constant curvature  $k=-1$. All necessary definitions from hyperbolic geometry can be found in the book \cite{Avs}.

By definition, a \emph{cyclic polygon} in the hyperbolic plane is a convex polygon inscribed in a circle, horocycle or one branch of equidistant curve. Useful information about cyclic polygons can be found in \cite{Walt1} and \cite{Walt2}. In particular, it is shown in \cite{Walt1} that any cyclic polygon in the hyperbolic plane is uniquely determined (up to isometry) by the ordered sequence of its side lengths. In addition, among all hyperbolic polygons with fixed positive side lengths there exist polygons of maximal area. Every such a maximal polygon is cyclic (see \cite{Walt2}).

The following four non-Euclidean versions of the Heron formula in the hyperbolic plane are known for a long time.

\begin{theorem}\label{ThGerLob} Let $T$ be a hyperbolic triangle with side lengths $a,b,c$. Then the area $S=S(T)$  is given by each of the following formulas
\begin{itemize}
\item[(i)]{Sine of 1/2 Area Formula}
$$\sin^2{\frac{S}{2}}=\frac{\sinh(s-a)\sinh(s-b)\sinh(s-c)\sinh(s)} {4\cosh^2{(\frac{a}{2})}\cosh^2{(\frac{b}{2})}\cosh^2{(\frac{c}{2})}};$$
\item[(ii)]{Tangent of 1/4 Area Formula}
$$\tan^2{\frac{S}{4}}=\tanh\left(\frac{s-a}{2}\right)\tanh\left(\frac{s-b}{2}\right)\tanh\left(\frac{s-c}{2}\right)\tanh\left(\frac{s}{2}\right);$$
\item[(iii)]{Sine of 1/4 Area Formula}
$$\sin^2{\frac{S}{4}}=\frac{\sinh(\frac{s-a}{2})\sinh(\frac{s-b}{2})\sinh(\frac{s-c}{2})\sinh(\frac{s}{2})} {\cosh{(\frac{a}{2})}\cosh{(\frac{b}{2})}\cosh{(\frac{c}{2})}};$$
\item[(iv)]{Bilinski Formula}
$$\cos{\frac{S}{2}}=\frac{\cosh a + \cosh b + \cosh c + 1} {4\cosh{(\frac{a}{2})}\cosh{(\frac{b}{2})}\cosh{(\frac{c}{2})}}.$$
\end{itemize}
\end{theorem}

The first two formulas are contained in the book \cite{Avs} (p.~66).
The third formula can be obtained by the squaring of the product of the first two. The forth one
was derived by Stanko Bilinski in \cite{Bilin} (see also \cite{Walt2}). 

It should be noted that the analogous formulas in spherical space are also known. For example, the spherical version of (i) is called the Cagnoli’s Theorem (see \cite{Todth}, sec.~100), (ii) is called the Lhuilier’s Theorem (see \cite{Todth}, sec.~102), (iii) and (iv) are proven in \cite{Todth} (sec.~103).

\subsection{Preliminary results for cyclic quadrilaterals}

We recall a few well known facts about  cyclic quadrilaterals. A convex Euclidean quadrilateral with interior angles  $A,B,C,D$ is cyclic if and only if $A+C=B+D=\pi$. A similar result for hyperbolic quadrilateral was obtained by V.~F.~Petrov \cite{Petr} and L.~Wimmer \cite{Wimmer}. They proved the following proposition.

\begin{proposition}\label{prop4-1}
A convex hyperbolic quadrilateral with interior angles $A,B,C,D$ is cyclic if and only if $A+C=B+D.$
\end{proposition}

\begin{figure}[h]
\begin{center}
\includegraphics[scale=0.4]{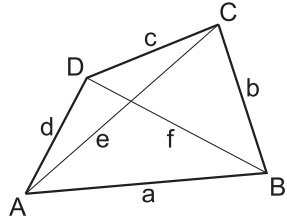}
\end{center}
\caption{Cyclic quadrilateral $Q$ with angles $A,B,C,D$, where $A+C=B+D$}
\end{figure}

Note that the sum of angles of a hyperbolic quadrilateral is less then $2\pi.$ Hence, for any cyclic  hyperbolic quadrilateral we have $A+C=B+D<\pi.$

It was shown in \cite{Walt1} that a cyclic $n$-gon is uniquely up to isometry determined by the lengths of its sides. Denote side and diagonal lengths of a quadrilateral as indicated on Fig. 7. Then, in the Euclidean case, a quadrilateral is cyclic if and only if \,$ef=a\,c+b\,d.$ This is the Ptolemy's theorem. A similar result for hyperbolic  quadrilateral is contained in the paper by J.~E.~Valentine \cite{Valentine}.

\begin{proposition}\label{prop2}
A convex hyperbolic quadrilateral with side lengths $a,b,c,d$ and diagonal lengths $e,f$ is cyclic if and only if
$$\sinh \frac{e}{2}\sinh \frac{f}{2}= \sinh \frac{a}{2}\sinh \frac{c}{2} + \sinh \frac{b}{2}\sinh \frac{d}{2}. $$
\end{proposition}

An important supplement to the Ptolemy's theorem is the following property of a cyclic quadrilateral in the Euclidean plane. Its side and diagonal lengths are related by the equation
\begin{equation}\label{Bret2}
\frac{e}{f}=\frac{a\,d+b\,c}{a\,b+c\,d}.
\end{equation}
Together with the Ptolemy's theorem this equation allows to express the diagonal lengths of a cyclic quadrilateral by its side lengths.

It was noted in \cite{GuoSonmez} that the above mentioned relationships between sides and diagonals of a cyclic quadrilateral are valid also in the hyperbolic geometry. To make them true, one can change the length side $a$ by the quantity $\displaystyle s(a)=\sinh \frac{a}{2}.$ In particular, formula (\ref{Bret2}) can be rewritten in the following way.

\begin{proposition}\label{prop4}
The side lengths $a,b,c,d$ and diagonal lengths $e,f$ of a cyclic hyperbolic quadrilateral are related by the following equation
$$\frac{s(e)}{s(f)}= \frac{s(a)s(d)+s(b)s(c)}{s(a)s(b)+s(c)s(d)}.$$
\end{proposition}

By making use of Propositions~\ref{prop2} and \ref{prop4} we derive the following formulas for the diagonal lengths $e,f$ of a cyclic hyperbolic quadrilateral. Then we have
\begin{equation}\label{eq4}
s^2(e)=\frac{s(a)s(d)+s(b)s(c)}{s(a)s(b)+s(c)s(d)}(s(a)s(c)+s(b)s(d)),
\end{equation}

\begin{equation}\label{eq5}
s^2(f)=\frac{s(a)s(b)+s(c)s(d)}{s(a)s(d)+s(b)s(c)}(s(a)s(c)+s(b)s(d)).
\end{equation}

Note that formulas (\ref{eq4}) and (\ref{eq5}) take a place also in the Euclidean and spherical geometries. In these cases, instead of function $s(a)$ one should take the functions $s(a)=a$ and $\displaystyle s(a)=\sin\frac{a}{2}$ respectively. See \cite{GuoSonmez} and \cite{DickSal} for the arguments in the spherical case.

All the above propositions will be used in the next section to obtain a few versions of the Brahmagupta formula for a cyclic hyperbolic quadrilateral.

\subsection{Area of cyclic hyperbolic quadrilateral}

In this section we consider the four versions of the Brahmagupta formula for a cyclic hyperbolic quadrilateral given by second named author \cite{Medn12}. They are generalizations of the respective statements (i)--(iv) of Theorem~\ref{ThGerLob}.

In particular, the first statement (i) has the following analog.

\begin{theorem}[Sine of $1/2$ area formula]\label{theo3}
Let $Q$ be a cyclic hyperbolic quadrilateral with side lengths $a,b,c,d$. Then the area $S=S(Q)$ is given by  the formula
$$
\sin^2\frac{S}{2}=\frac{\sinh(s-a)\sinh(s-b)\sinh(s-c)\sinh(s-d)}
{4\cosh^2\frac{a}{2}\cosh^2\frac{b}{2}\cosh^2\frac{c}{2}\cosh^2\frac{d}{2}}(1-\varepsilon),
$$
where \,$\varepsilon=\displaystyle{\frac{\sinh\frac{a}{2}\sinh\frac{b}{2}\sinh\frac{c}{2}\sinh\frac{d}{2}}
{\cosh\frac{s-a}{2}\cosh\frac{s-b}{2}\cosh\frac{s-c}{2}\cosh\frac{s-d}{2}}}$ and \,$\displaystyle{s=\frac{a+b+c+d}{2}}$.
\end{theorem}

Note that the number $\varepsilon$ vanishes if $d=0$. In this case, we get formula (i) again.

The second statement (ii) for the case of a hyperbolic quadrilateral has the following form.

\begin{theorem}[Tangent of $1/4$ area formula]\label{theo2}
Let $Q$ be a cyclic hyperbolic quadrilateral with side lengths $a,b,c,d$. Then the area $S=S(Q)$ is given by  the formula
$$
\tan^2\frac{S}{4}=\frac{1}{1-\varepsilon}\tanh\frac{s-a}{2}\tanh\frac{s-b}{2}\tanh\frac{s-c}{2}\tanh\frac{s-d}{2},
$$
where $s$ and $\varepsilon$ are the same as in Theorem~\ref{theo3}.
\end{theorem}

It follows from Theorem~\ref{theo3} that for any $a,b,c,d \neq 0$ we have $1-\varepsilon>0$ and $\varepsilon>0$. Hence,  $0<\varepsilon<1$. Taking into account these inequalities as an immediate consequence of theorems~\ref{theo3} and \ref{theo2} we obtain the following corollary.

\begin{corollary} For any cyclic hyperbolic quadrilateral the following inequalities take a place
$$\sin^2{\frac{S}{2}}<\frac{\sinh(s-a)\sinh(s-b)\sinh(s-c)\sinh(s-d)} {4\cosh^2{(\frac{a}{2})}\cosh^2{(\frac{b}{2})}\cosh^2{(\frac{c}{2})}\cosh^2{(\frac{d}{2})}}$$
and
$$\tan^2{\frac{S}{4}}>\tanh\left(\frac{s-a}{2}\right)\tanh\left(\frac{s-b}{2}\right)\tanh\left(\frac{s-c}{2}\right) \tanh\left(\frac{s-d}{2}\right).$$
\end{corollary}

By squaring the product of formulas in Theorems~\ref{theo3} and \ref{theo2} we obtain the following result. It can be considered as a direct generalization of the statement (iii) in Theorem~\ref{ThGerLob}.

\begin{theorem}[Sine of $1/4$ area formula]\label{theo4}
Let $Q$ be a cyclic hyperbolic quadrilateral with side lengths $a,b,c,d$. Then the area $S=S(Q)$ is given by  the formula
$$
\sin^2\frac{S}{4}=\displaystyle{\frac{\sinh\frac{s-a}{2}\sinh\frac{s-b}{2}\sinh\frac{s-c}{2}\sinh\frac{s-d}{2}}
{\cosh\frac{a}{2}\cosh\frac{b}{2}\cosh\frac{c}{2}\cosh\frac{d}{2}}},
$$
where $\displaystyle{s=\frac{a+b+c+d}{2}}$.
\end{theorem}

The analogous formula in spherical space can be obtained by replacing $\sinh$ and $\cosh$ with $\sin$ and $\cos$ correspondingly (see \cite{MCel}, p.~182, proposition~5).

Consider a circumscribed quadrilateral $Q$ with side lengths $a,b,c,d$ (Fig.~8). In this case we have
$s-a=c,\, s-b=d,\, s-c=a,\, s-d=b$. As a result we obtain the following assertion.

\bigskip
\begin{figure}[h]\label{main}
\begin{center}
\includegraphics[scale=0.4]{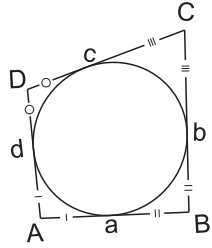}
\end{center}
\caption{Circumscribed cyclic quadrilateral $Q$ with side lengths $a,b,c,d$}
\end{figure}

\begin{corollary}[Brahmagupta formula for bicentric quadrilateral]
Let $Q$ be a bicentric (i.e. inscribed and circumscribed) hyperbolic quadrilateral with side lengths $a,b,c,d$. Then the area $S=S(Q)$ is given by the formula
$$\sin^2{\frac{S}{4}}=\tanh\left(\frac{a}{2}\right)\tanh\left(\frac{b}{2}\right)\tanh\left(\frac{c}{2}\right)\tanh\left(\frac{d}{2}\right).$$
\end{corollary}

The analogous formula in spherical space can be given by replacing $\tanh$ with $\tan$ (see \cite{MCel}, p. 46). An Euclidean version of this result is well known (see, e.g., \cite{Ivan}). In this case $$S^2=a\,b\,c\,d.$$

The next theorem \cite{Medn12} presents a version of the Bilinski formula for a cyclic quadrilateral.

\begin{theorem}[Bilinski formula]\label{theo1}
Let $Q$ be a cyclic hyperbolic quadrilateral with side lengths $a,b,c,d$. Then the area $S=S(Q)$ is given by the formula
$$\cos{\frac{S}{2}}=\frac{\cosh a + \cosh b + \cosh c  + \cosh d -{4\sinh{(\frac{a}{2})}\sinh{(\frac{b}{2})}\sinh{(\frac{c}{2})}}\sinh{(\frac{d}{2})}} {4\cosh{(\frac{a}{2})}\cosh{(\frac{b}{2})}\cosh{(\frac{c}{2})}\cosh{(\frac{d}{2})}}.$$
\end{theorem}


\subsection{Proof of Brahmagupta formulas}

Consider a cyclic hyperbolic quadrilateral $Q$ with side lengths $a,b,c,d$ and interior angles $A,B,C,D$ shown on Fig.~7. By the Gauss--Bonnet formula we get the area $$S=S(Q)=2\pi-A-B-C-D.$$

To prove Theorem~\ref{theo3} let us find the quantities $\displaystyle\sin^2\frac{S}{4}$ and $\displaystyle\cos^2\frac{S}{4}$ in terms of $a,b,c,d$. Since $A+C=B+D$ (see Proposition~\ref{prop4-1}) we have
$$2\sin^2\frac{S}{4}=1-\cos\frac{S}{2}=1-\cos(\pi-(A+C))=1+\cos(A+C).$$

Hence,
\begin{equation}
\label{eq13}
\sin^2\frac{S}{4}=\frac{1+\cos A\,\cos C-\sin A\,\sin C}{2}.
\end{equation}

Now we show that $\cos A,\cos C$ and the product $\sin A\cdot\sin C$ can be expressed in terms of elementary functions in $a,b,c,d$. To find $\cos A$ we use the Cosine rule for hyperbolic triangle $ABD$.
\begin{equation}
\label{eq3}
\cos A=\frac{\cosh a\,\cosh d-\cosh f}{\sinh a\,\sinh d}.
\end{equation}

We note that $\cosh f =2\,s^2(f)+1, \cosh a=2\,s^2(a)+1$ and $\cosh d=2\,s^2(d)+1.$ Putting these identities into equations (\ref{eq5}) and (\ref{eq3}) we express $\cos A$ in terms of $a,b,c,d.$ After straightforward calculations  we obtain
\begin{equation}
\label{eq6}
\cos A=\frac{s^2(a)-s^2(b)-s^2(c)+s^2(d)+2s(a)\,s(b)\,s(c)\,s(d)+2s^2(a)\,s^2(d)}{2\,(s(a)\,s(d)+s(b)\,s(c))\cosh\frac{a}{2}\cosh\frac{d}{2}}.
\end{equation}

In a similar way we get the formula
\begin{equation}
\label{eq66}
\cos C=\frac{-s^2(a)+s^2(b)+s^2(c)-s^2(d)+2s(a)\,s(b)\,s(c)\,s(d)+2s^2(b)\,s^2(c)}
{2\,(s(a)\,s(d)+s(b)\,s(c))\cosh\frac{b}{2}\cosh\frac{c}{2}}. \quad
\end{equation}

Note that $\sin A\,\sin C > 0$ and $\sin^2 A\,\sin^2 C= (1-\cos^2 A)(1-\cos^2 C)$. Then we take a square root in the latter equation, where $\cos A$ and $\cos C$ are given by formulas  (\ref{eq6}) and (\ref{eq66}). Thus we express a product $\sin A\cdot\sin C$ in terms of $a,b,c,d$. Substituting $\cos A, \cos C$ and $\sin A\cdot\sin C$ into (\ref{eq13}) and simplifying we get

\begin{equation}
\label{eq14}
\sin^2\frac{S}{4}=\frac{\sinh\frac{-a+b+c+d}{4}\sinh\frac{a-b+c+d}{4}\sinh\frac{a+b-c+d}{4}\sinh\frac{a+b+c-d}{4}}
{\cosh\frac{a}{2}\cosh\frac{b}{2}\cosh\frac{c}{2}\cosh\frac{d}{2}}.
\end{equation}

This proves Theorem~\ref{theo4}.
\medskip

In a similar way, from identity $\displaystyle2\cos^2\frac{S}{4}=1+\cos\frac{S}{2}=1-\cos(A+C)$ we have
\begin{equation}
\label{eq15}
\cos^2\frac{S}{4}=\frac{\cosh\frac{a+b-c-d}{4}\cosh\frac{a-b+c-d}{4}\cosh\frac{a-b-c+d}{4}\cosh\frac{a+b+c+d}{4}}
{\cosh\frac{a}{2}\cosh\frac{b}{2}\cosh\frac{c}{2}\cosh\frac{d}{2}}.\end{equation}

The following lemma can be easy proved by straightforward calculations.

\begin{lemma}\label{lem1} 
The expression
$$
H=\frac{\cosh\frac{a+b-c-d}{4}\cosh\frac{a-b+c-d}{4}\cosh\frac{a-b-c+d}{4}\cosh\frac{a+b+c+d}{4}
}{\cosh\frac{-a+b+c+d}{4}\cosh\frac{a-b+c+d}{4}\cosh\frac{a+b-c+d}{4}\cosh\frac{a+b+c-d}{4}}
$$
can be rewritten in the form $H=1-\varepsilon,$ \textit{where} $$\varepsilon=\frac{\sinh\frac{a}{2}\sinh\frac{b}{2}\sinh\frac{c}{2}\sinh\frac{d}{2}}
{\cosh\frac{s-a}{2}\cosh\frac{s-b}{2}\cosh\frac{s-c}{2}\cosh\frac{s-d}{2}}\textit{ and }s=\frac{a+b+c+d}{2}.$$
\end{lemma}

\bigskip
Taking four times product of equations (\ref{eq14}) and (\ref{eq15}) we have

\begin{equation}
\label{eq16}
\sin^2\frac{S}{2}=\frac{\sinh\frac{-a+b+c+d}{2}\sinh\frac{a-b+c+d}{2}\sinh\frac{a+b-c+d}{2}\sinh\frac{a+b+c-d}{2}}
{4\cosh^2\frac{a}{2}\cosh^2\frac{b}{2}\cosh^2\frac{c}{2}\cosh^2\frac{d}{2}}\cdot H,
\end{equation}
where $H$ is the same as in Lemma~\ref{lem1}.

Then the statement of Theorem~\ref{theo3} follows from equation (\ref{eq16}), Lemma~\ref{lem1} and the evident identity $\displaystyle s-a=\frac{-a+b+c+d}{2}.$

\bigskip
To prove Theorem~\ref{theo2} we divide (\ref{eq14}) by (\ref{eq15}). As a result we have

\begin{equation}
\label{eq17}
\tan^2\frac{S}{4}=\frac{\sinh\frac{-a+b+c+d}{4}\sinh\frac{a-b+c+d}{4}\sinh\frac{a+b-c+d}{4}\sinh\frac{a+b+c-d}{4}}
{\cosh\frac{a+b-c-d}{4}\cosh\frac{a-b+c-d}{4}\cosh\frac{a-b-c+d}{4}\cosh\frac{a+b+c+d}{4}}.
\end{equation}

Hence, applying Lemma~\ref{lem1} we obtain the statement of Theorem~\ref{theo2}.

\bigskip
Finally, the Bilinski formula (Theorem~\ref{theo1}) follows from the identity $\displaystyle\cos{\frac{S}{2}}=\cos^2{\frac{S}{4}}-\sin^2{\frac{S}{4}}$ and the above mentioned equations (\ref{eq14}) and (\ref{eq15}).


\subsection{Area of hyperbolic trapezoid}

In this section we give a formula for the area of a hyperbolic trapezoid in terms of its side lengths. 

A convex hyperbolic quadrilateral with interior angles $A,B,C,D$ is called a \emph{trapezoid} if $A+B=C+D$ (see Fig.~9). This definition is also valid for Euclidean case. 

\bigskip
\begin{figure}[h]\label{main}
\begin{center}
\includegraphics[scale=0.4]{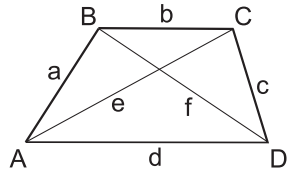}
\end{center}
\caption{Trapezoid $T=T(a,b,c,d)$}
\end{figure}

Denote the lengths of sides and diagonals as shown on Fig.~9. We assume that $b\ne d$. Otherwise, in the case $b=d$, the area of trapezoid $T$ is not determined by lengths of its sides $a,b,c,d$. The next formula was obtained by Dasha Sokolova \cite{Dasha}.

\begin{theorem}[Sokolova, 2012]
Let $T$ be a hyperbolic trapezoid with side lengths $a,b,c,d$. Then the area $S=S(T)$ is given by the formula
$$
\tan^2\frac{S}{4}=
\frac{\sinh^2\frac{b+d}{2}\sinh\frac{a+b-c-d}{4}\sinh\frac{a+b+c-d}{4}\sinh\frac{-a+b+c-d}{4}\sinh\frac{a-b+c+d}{4}}
{\sinh^2\frac{b-d}{2}\cosh\frac{a-b-c-d}{4}\cosh\frac{a-b+c-d}{4}\cosh\frac{a+b-c+d}{4}\cosh\frac{a+b+c+d}{4}}.
$$
\end{theorem}

\begin{remark}
The following formula gives the area of an Euclidean trapezoid in terms of its side lengths.
$$
{S_E}^2=\frac{(b+d)^2(a+b-c-d)(a+b+c-d)(-a+b+c-d)(a-b+c+d)}{16\,(b-d)^2}.
$$

Note that $\displaystyle\tan^2\frac{S}{4}\sim\left(\frac{S_E}{4}\right)^2$ for sufficiently small values of $a,b,c,d$.
\end{remark}

\section{Acknowledgments}

The authors would like to express their gratitude to a referee for careful consideration of the paper and useful remarks and suggestions.

\end{document}